\input amstex
\documentstyle {amsppt}
\UseAMSsymbols \vsize 18cm \widestnumber\key{ZZZZZ}

\catcode`\@=11
\def\displaylinesno #1{\displ@y\halign{
\hbox to\displaywidth{$\@lign\hfil\displaystyle##\hfil$}&
\llap{$##$}\crcr#1\crcr}}
\def\ldisplaylinesno #1{\displ@y\halign{
\hbox to\displaywidth{$\@lign\hfil\displaystyle##\hfil$}&
\kern-\displaywidth\rlap{$##$} \tabskip\displaywidth\crcr#1\crcr}}
\catcode`\@=12

\refstyle{A}

\let \ol=\overline

\let \ti=\widetilde

\font\sit=cmti10 at 7.5pt

 \font\srm=cmr10 at 7.5pt

\font \fin=lasy8 at 15.4 pt

\def \s{\mathop{\hbox {\rm s}}\nolimits}

\def \GL{\mathop{\hbox{\rm GL}}\nolimits}

\def \Ad{\mathop{\hbox{\rm Ad}}\nolimits}
\def \rk{\mathop{\hbox{\rm rk}}\nolimits}

\def \GL{\mathop{\hbox{\rm GL}}\nolimits}
\def \SL{\mathop{\hbox{\rm SL}}\nolimits}

\topmatter
\title Unipotent orbits and local $L$-functions \endtitle

\rightheadtext{Unipotent orbits and local $L$-functions}
\author Volker Heiermann \endauthor

\address Institut f\"ur Mathematik, Humboldt-Universit\"at zu
Berlin, Unter den Linden 6, 10099 Berlin, Allemagne \endaddress

\email heierman\@mathematik.hu-berlin.de \endemail

\abstract In a previous article (\it Orbites unipotentes et
p\^oles d'ordre maximal de la fonction $\mu $ de Harish-Chandra,
\rm to appear in Canad. J. Math.), we have assumed the existence
of the local Langlands correspondence for supercuspidal
representations and deduced from this a local Langlands
correspondence for discrete series representations and beyond
(without going into the structure of the $L$-packets). The aim of
the present article is to show that this extension of the local
Langlands correspondence for supercuspidal representations (and
some of the assumptions in the article above) is compatible with
the theory of $L$-functions due to Langlands-Shahidi.
 \rm

\endabstract

\endtopmatter
\document

Let $G$ be the group of $F$-points of a connected reductive group
defined over a non archimedean local field. In \cite{H2} we have
assumed the local Langlands correspondence for supercuspidal
representations of $F$-Levi-subgroups of $G$ and deduced from this
a local Langlands correspondence for discrete series
representations of $G$ and beyond (without going into the
structure of the $L$-packets).  The aim of this note is to show
that the results and assumptions in \cite{H2} are compatible with
the theory of $L$-functions of Langlands-Shahidi. This theory
applies at this moment to generic representations of $F$-points of
quasi-split connected reductive groups. It has been established
until now only for $F$ of characteristic $0$. So we have to make
this assumption, too, and suppose in the sequel that $G$ is
quasi-split.

Let us be more precise. Let $P=MU$ be an $F$-parabolic subgroup of
$G$. Denote by $\Sigma _{red}(P)$ the set of reduced roots in
$Lie(U)$ of the maximal split torus $A_M$ in the center of $M$.
Recall that to each $\alpha\in\Sigma _{red}(P)$ corresponds a
semi-standard $F$-Levi subgroup $M_{\alpha }$ of $G$, which
contains $M$ as maximal Levi subgroup. One identifies $\Sigma
_{red}(P)$ to a subset of the dual $a_M^*$ of the real
 Lie algebra of $A_M$. There is a natural way
to attach to an element $\lambda $ of the complexification of
$a_M^*$ a character $\chi _{\lambda }$ of $M$ \cite{H2, 0.6}. If
$\lambda =s\alpha $, $s\in\Bbb C$, and $m\in A_M$, one has $\chi
_{s\alpha }(m)=\vert\alpha(m)\vert_F^s$, where $\vert\cdot\vert
_F$ denotes the normalized absolute value of $F$.

Let $\sigma $ be an irreducible unitary supercuspidal generic
representation of $M$ and $W_F$ the Weil group of $F$. In
\cite{H2} we have assumed that one can attach to $\sigma $ an
admissible homomorphism $\psi _{\sigma }:W_F\times\SL _2(\Bbb
C)\rightarrow\ ^LM$ (see \cite{H2} for the precise definition of
the Langlands $L$-group and an admissible homomorphism used here),
verifying some properties, coming from the conjectural local
Langlands correspondence. As in particular it is believed that
${\psi _{\sigma }}_{\vert \SL _2(\Bbb C)}$ is trivial, when $\psi
_{\sigma }$ is attached to a generic supercuspidal representation,
we will assume this here, too. The assumption \cite{H2,4.3}
simplifies then considerably and reads (with $q$ the number of
elements in the residue field of $F$, referring to \cite{H2, 3.5}
for the notion of "q-distinguished")

\null (LM) \it For each root $\alpha\in\Sigma _{red}(P)$, the
Harish-Chandra $\mu $-function $s\mapsto\mu^{M_{\alpha
}}(\sigma\otimes\chi_{s\alpha })$ (see \cite{W} for the definition
of this function) has a pole in a real number $s_0>0$, if and only
if $\alpha (q)^{s_0}$ is $q$-distinguished in the connected
centralizer of the image of $\psi _{\sigma }$ and this group is
not a torus. \rm

\null Fix a non trivial additive character $\psi _F$ of $F$. In
\cite{Sh} Shahidi (proving a conjecture of Langlands) has
associated to an irreducible smooth generic representation $\sigma
$ of $M$ a set of complex functions $\{s\mapsto\gamma (s,\sigma
,r_i,\psi _F), 1\leq i\leq m\}$. From them he deduces canonically
$L$-functions $L(\s,\sigma ,r_i)$ and $\epsilon $-factors
$\epsilon(s,\sigma ,r_i,\psi _F)$ (see also {\bf 1.3} for more
details). As the maps $r_i\circ\psi _{\sigma }$ are
representations of the Weil-Deligne group, the Artin $L$-functions
$L(s,r_i\circ\psi _{\sigma })$ and $\epsilon $-functions $\epsilon
(s,r_i\circ\psi _{\sigma },\psi _F)$ are defined and one derives
from them $\gamma (s,r_i\circ\psi _{\sigma },\psi _F)$ as above
(see {\bf 1.4 - 1.5} for more details).

\null Our first result is, that the assumption $(LM)$ is
satisfied, if $\sigma $ and $\psi _{\sigma }$ have the same
$L$-functions w.r.t. each $M_{\alpha }$, $\alpha\in \Sigma
_{red}(P)$. We get also a converse under some condition on the
$L$-functions attached to $\psi _{\sigma }$.

\null Under the assumption (LM) we have in \cite{H2} associated to
each elliptic admissible homomorphism $\psi :W_F\times\SL _2(\Bbb
C)\rightarrow\ ^LG$ an irreducible square-integrable
representation $\pi $ of $G$, and vice-versa. Our next result is
that $\psi $ and $\pi $ have same $\gamma $-functions if they
correspond to each other by this correspondence. We show also that
this property remains true, if one extends the correspondence to
arbitrary admissible homomorphisms $\psi $ and arbitrary smooth
irreducible representations $\pi $ of $G$, as done in the last
section of \cite{H2}.

\null We finish by a discussion of the general case of non generic
representations and non quasi-split groups, in taking into account
the conjectural framework in \cite{Sh, 9.}.

\null We refer to the introduction of \cite{H2} for information of
the actual state of the local Langlands conjectures.

\null {\bf 1. Notations and preliminaries:}

\null{\bf 1.1.} We denote by $I_F$ the inertial subgroup of $W_F$,
by $Fr$ a geometric Frobenius automorphism of $F$ \cite{De} and
normalize the reciprocity map in local class field theory so that
$\vert Fr\vert _F=q^{-1}$.

To simplify the notations, we will denote by $\Im (s)$ the
imaginary part of a complex number $s$ multiplied by $\sqrt{-1}$.

\null {\bf 1.2} We fix a minimal $F$-parabolic subgroup
$P_0=M_0U_0$ of $G$ and a maximal $F$-split torus $A_0$ contained
in $M_0$. We denote by $\Sigma $ the set of roots of $A_0$ in
$Lie(G)$ and by $\Delta $ the set of simple roots with respect to
$P_0$. If $P=MU$ is a standard parabolic of $G$ (i.e. $P\supseteq
P_0$), $\alpha\in\Sigma _{red}(P)$, we note $P_{\alpha }$ the
standard parabolic $P\cap M_{\alpha }$ of $M_{\alpha }$ and
$U_{\alpha }=U\cap M_{\alpha }$.

\null{\bf 1.3} Let $P=MU$ be a maximal standard $F$-parabolic
subgroup of $G$, $\rho $ half of the sum of the roots in $\Sigma $
whose root space spans $Lie(U)$ and $\alpha $ the unique root in
$\Delta $ which does not lie in the root subsystem of $\Sigma $
corresponding to $M$. Put $\ti {\alpha }=\langle\rho ,\alpha^{\vee
}\rangle^{-1}\rho $.

Denote by $r$ the adjoint action of $^LM$ on $Lie(^LU)$ and
$$V_i=\{X_{\beta ^{\vee }}\in Lie(^LU)\vert\ \langle\ti{\alpha }
,\beta ^{\vee }\rangle=i\}.$$ (Here $Lie(^LU)$ has been decomposed
into weight spaces relative to the roots with respect to the
action of the connected center of $^LT$, which equals $^LA_0$.)
The spaces $V_i$ are invariant by $r$. Denote by $r_i$ the
restriction of $r$ to $V_i$. One has a decomposition $r=\oplus
_{i=1}^mr_i$ with some integer $m\geq 1$, called the length of
$r$. The components $r_i$, $1\leq i\leq m$, are irreducible
\cite{Sh}.

Let $\sigma $ be a smooth irreducible generic representation of
$M$. Fix a non trivial additive character $\psi _F$ of $F$. In
\cite{Sh} Shahidi (proving a conjecture of Langlands) has
associated to $\sigma $ a set of complex functions
$\{s\mapsto\gamma (s,\sigma ,r_i,\psi _F), 1\leq i\leq m\}$. If
$\sigma $ is tempered, he deduces from them canonically
$L$-functions $L(s,\sigma ,r_i)$ and $\epsilon $-factors
$\epsilon(s,\sigma ,r_i,\psi _F)$ in the following way: Denote by
$P_{\sigma ,i}$ the unique polynomial satisfying $P_{\sigma
,i}(0)=1$ such that $P_{\sigma ,i}(q^{-s})$ is the numerator of
$\gamma (s,\sigma ,r_i,\psi _F)$ (in particular $P_{\sigma
,i}(q^{-s})$ has the same zeros than $\gamma (s,\sigma ,r_i,\psi
_F)$). Then
$$L(s,\sigma ,r_i):=P_{\sigma ,i}(q^{-s})^{-1},\qquad
L(s,\sigma ,\ti{r}_i):=P_{\ti{\sigma },i}(q^{-s})^{-1}$$
$$\hbox{\rm and}\qquad\epsilon(s,\sigma ,r_i,\psi_F):=
(L(s,\sigma ,r_i)/L(1-s,\sigma ,\ti{r}_i))\gamma (s,\sigma
,r_i,\psi _F)$$ for $1\leq i\leq m$, where $\ti{\sigma }$ and
$\ti{r}_i$ are the contragredient representations.

As $\gamma (s,\sigma\otimes\chi _{s'\ti{\alpha }},r_i,\psi
_F)=\gamma (s+s',\sigma ,r_i,\psi _F)$ by \cite{Sh, (3.12)},
$L(s,\sigma ,r_i)$ and $\epsilon(s,\sigma ,r_i,\psi_F)$ are also
defined, if $\sigma $ is only quasi-tempered.

\null The following properties hold:

(1.3.1) $L(s,\sigma ,r_i)=1$ for $3\leq i\leq m$, if $\sigma $ is
supercuspidal \cite{Sh, 7.5};

(1.3.2) Suppose that $P$ is associated to its opposite parabolic
subgroup $\ol{P}$ and that $\sigma $ is unitary and supercuspidal.
(We will later say that $P$ is self-conjugated.) Denote by $w$ a
representative of an element of the Weyl group that conjugates $P$
and $\ol{P}$. Then the Harish-Chandra $\mu $-function (see
\cite{W} for the definition of this function) verifies (with $\sim
$ meaning equality up to a monomial in $q^{-s}$)
$$\eqalign{\mu(\sigma\otimes\chi _{s\ti{\alpha }})&\sim
{P_{\sigma ,1}(q^{-s})P_{\sigma ,2}(q^{-2s})P_{w\sigma
,1}(q^{s})P_{w\sigma ,2}(q^{2s})\over P_{\ti{\sigma }
,1}(q^{-(1-s)})P_{\ti{\sigma },2}(q^{-(1-2s)})P_{w\ti{\sigma }
,1}(q^{-(1+s)})P_{w\ti{\sigma },2}(q^{-(1+2s)})}\cr &=
{L(1-s,\sigma ,\ti{r_1})L(1-2s,\sigma ,\ti{r_2})L(1+s,w\sigma
,\ti{r_1})L(1+2s,w\sigma ,\ti{r_2}) \over L(s,\sigma ,r_1)
L(2s,\sigma ,r_2)L(-s,w\sigma ,r_1)L(-2s,w\sigma ,r_2)}\cr}$$
\cite{Sh, 1.4 and 7.6}.

\null (1.3.3) $$\ol{L(s,\sigma ,r_i)}=L(\ol{s},\sigma ,\ti{r_i})=
L(\ol{s},\ti{\sigma },r_i)\eqno{\hbox{\rm \cite{Sh, 7.8 and p.
308}.}}$$

(1.3.4) $$L(s,\sigma\otimes\chi _{s'\ti{\alpha }},r_i)=L
(s+s',\sigma ,r_i)\eqno{\hbox{\rm \cite{Sh, (3.12)}.}}$$

\null (1.3.5) If $\sigma $ is unitary and supercuspidal, all poles
of $L(\cdot ,\sigma ,r_i)$ have real part $0$ \cite{Sh, 7.3}.

\null (1.3.6) If $\sigma $ is supercuspidal, the poles of $L(\cdot
 ,\sigma ,r_i)$ are simple. (This is because of (1.3.2) and the
simplicity of the poles of Harish-Chandra's $\mu $-function
\cite{H1, remark in the proof of 4.1}.)

\null If $(r',V')$ is a sub-representation of $r$, one defines
$\gamma (\cdot ,\sigma ,r',\psi _F)=\prod _{i,V_i\subseteq
V'}\gamma (\cdot ,\sigma ,$ $r_i,\psi _F)$ and in the same way
$L(\cdot ,\sigma ,r')$ and $\epsilon (\cdot ,\sigma ,r',\psi _F) $

If $\pi $ is a general generic smooth irreducible representation
of $M$, then the $L$-functions $L(\cdot ,\pi ,r_i)$ are defined in
the following way \cite{Sh, p. 308}: by Langlands' classification
there is a standard $F$-parabolic subgroup $P_1=M_1U_1$ of $G$
with $M_1\subseteq M$ and an irreducible quasi-tempered
representation $\tau $ of $M_1$, such that $\pi $ is the unique
sub-representation of $i_{P_1\cap M}^M\tau $. By \cite{R, Theorem
2}, the quasi-tempered representation $\tau $ is generic. Denote
by $\kappa _1$ the inclusion $^LM_1\rightarrow\ ^LM$, by $r_{1,i}$
the composition $r\circ\kappa _1$ and, for $\alpha\in\Sigma
_{red}(P_1)-\Sigma _{red}(P_1\cap M)$, by $r_{1,i,\alpha }$ the
restriction $^LM_1\rightarrow Lie(^LU_{1,\alpha })$ of $r_{1,i}$.
The $L$-functions $L(\cdot ,\tau ,r_{1,i,\alpha })$ w.r.t.
$M_{1,\alpha }$ are defined by analytic continuation from the
tempered case.  The $L$-function associated to $\pi $ and $r_i$ is
$$L(\cdot, \pi ,r_i)=\prod _{\alpha\in\Sigma
_{red}(P_1)-\Sigma _{red}(P_1\cap M)} L(\cdot, \tau ,r_{1,i,\alpha
}).\leqno{(1.3.7)}$$ The corresponding $\epsilon $-factor is
deduced from $L(\cdot, \pi ,r_i)$ and $\gamma (\cdot, \pi
,r_i,\psi _F)$ by the same equation as in the tempered case.

Consider finally an arbitrary standard parabolic subgroup $P=MU$
of $G$. Denote still by $\rho $ half of the sum of the roots in
$\Sigma $ that generate $U$. For each $\beta $ with $X_{\beta
^{\vee }}\in Lie(^LU)$, $\langle2\rho,\beta ^{\vee }\rangle $ is a
positive integer. Let $1\leq a_1<a_2<\dots <a_m$ be the different
values. Following \cite{Sh}, we define
$$V_i:=\{X_{\beta ^{\vee }}\in Lie(^LU)\vert \langle 2\rho,\beta
^{\vee }\rangle =a_i\}$$ and denote by $r_i$ the restriction of
the adjoint representation $r:\ ^LM\rightarrow Lie(^LU)$ to $V_i$.
This definition agrees with the one above for $P$ maximal. For
$\alpha\in\Sigma _{red}(P)$, let $r_{i,\alpha }:\ ^LM\rightarrow
Lie(^LU_{\alpha })$ be the restriction of $r_i$. If $r_{\alpha }$
denotes the adjoint representation $^LM\rightarrow Lie(^LU_{\alpha
})$, then it follows from elementary properties of root systems
that $r_{i,\alpha }=r_{\alpha ,i}$ (with $r_{\alpha ,i}$ defined
relative to the maximal parabolic subgroup $P\cap M_{\alpha }$ of
$M_{\alpha }$ as above). Let $\pi $ be a general generic
irreducible smooth representation of $M$. For $\alpha\in\Sigma
_{red}(P)$, denote by $\gamma(\cdot ,\pi ,r_{i,\alpha },\psi _F)$
the $\gamma $-function of $\pi $ defined relative to $M_{\alpha }$
and $P\cap M_{\alpha }$. Then, by definition,
$$\gamma (\cdot ,\pi ,r_i,\psi _F)=\prod _{\alpha\in\Sigma
_{red}(P)}\gamma (\cdot,\pi ,r_{i,\alpha },\psi _F)$$ \cite{Sh, p.
307, l.15 -20}. The $L$- and $\epsilon $-factors of $\pi $
relative to $P$ are defined in the same way as product of $L$- and
$\epsilon $-factors attached to $\alpha\in\Sigma _{red}(P)$.

If $r'$ is an arbitrary sub-representation of $r$, then one
defines local factors for $r'$ in the same way than for maximal
$P$.

\null {\bf 1.4} Recall the definition of the Artin $L$-function
\cite{De}. An admissible homomorphism $\psi :W_F\times\SL_2(\Bbb
C)\rightarrow \GL_n(\Bbb C)$ can be written as direct sum of
twists of elliptic admissible homomorphisms. As the Artin
$L$-functions are additive and behave well under unramified twists
(i.e. $L(s+s',\psi )= L(s,\psi\ \vert\det\vert _F^{s'})$), it is
enough to give the definition for $\psi $ elliptic. Let $N$ be the
nilpotent $n\times n$-matrix, such that $\psi (\pmatrix 1 & 1
\\ 0 & 1\endpmatrix)=\exp (N)$. Identify $N$ with the
corresponding nilpotent endomorpism  of $V:=\Bbb C^n$. As $\psi $
is elliptic, the restriction $\psi _0$ of $\psi $ to $W_F$ is a
multiple of an irreducible representation and the subspace
$\ker(N)$ is an irreducible component. If $m$ is the multiplicity
of $\psi _0$ in $\psi $, one has
$$L(s,\psi )=\det (I-\psi (Fr)q^{-s-m+1}\vert \ker(N)^{I_F})^{-1}.$$
Remark that the action of $\psi _0$ on $V^{I_F}$ is an unramified
character. As $\ker(N)$ is an irreducible component of $\psi _0$,
either the representation $\psi _0$ is itself an unramified
character or $V^{I_F}=0$. So $L(s,\psi )=1$, if $\psi $ is
ramified. Otherwise $\dim(\ker(N)^{I_F})=1$. So, if $q^{s_0}$ is
the proper value for the action of $\psi (Fr)$ on $V^{I_F}$, then
we have $L(s,\psi )=(1-q^{s_0-n+1-s})^{-1}$ in this case. Remark
that $\Re (s_0)=0$, if $\psi (W_F)$ is relatively compact.

\null {\bf 1.5} The $\gamma $-, $L$- and $\epsilon $-factors
should be preserved by the (in general) conjectural local
Langlands correspondence. More precisely, let $\psi _{\sigma
}:W_F\times \SL_2(\Bbb C)\rightarrow\ ^LM$ be the conjectural
admissible homomorphism attached to the generic irreducible smooth
 representation $\sigma $. (It is in particular assumed that $\psi
_{\sigma }(W_F)$ is relatively compact, when $\sigma $ is
tempered.)

Then we should have
$$L(s,\sigma ,r_i)=L(s,r_i\circ\psi _{\sigma }) \qquad \hbox{\rm
and}\qquad \epsilon (s,\sigma ,r_i,\psi _F)=\epsilon
(s,r_i\circ\psi _{\sigma },\psi _F).$$ Here $L$- and $\epsilon
$-factors on the Galois side are the Artin $L$- and
$\epsilon$-functions defined by Deligne \cite{De}. If one defines
$\gamma (s,r_i\circ\psi _{\sigma }, \psi _F)$ by the corresponding
equation on the Galois side, one gets also
$$\gamma (s,\sigma ,r_i,\psi _F)=\gamma (s,r_i\circ\psi
_{\sigma },\psi _F).$$ Remark that, as $\epsilon (s,r_i\circ\psi
_{\sigma },\psi _F)$ is a monomial in $q^{-s}$ \cite{De}, $L(s,r
\circ\psi _{\sigma })^{-1}$ is the unique polynomial in
$z=q^{-s}$, which takes value $1$ in $z=0$ and which is the
numerator of $\gamma (s,r_i\circ\psi _{\sigma },\psi _F)$.

So, in particular, if $\sigma $ is tempered, the equality of
$\gamma $-factors implies the equality of $L$- and $\epsilon
$-factors.

\null \null{\bf 2.} We will now start to prove that in the generic
case the assumption $(LM)$ in \cite{H2} is implied by an equality
of $L$-functions (referring to \cite{H2, 3.5} for the notion of
"q-distinguished"), establishing also a kind of converse.

The lemma below is a reformulation of results in \cite{Sh}.

\null{\bf 2.1 Lemma:} \it Let $P=MU$ be a maximal standard
$F$-parabolic subgroup of $G$ and let $\sigma $ be a unitary
irreducible generic supercuspidal representation of $M$. Then, for
any $s\in\Bbb C$, $\mu (\sigma\otimes\cdot )$ has a pole in $\chi
_{s\ti{\alpha }}$, if and only if $\Im(s)$ is a pole of $L(\cdot
,\sigma ,r_i)$ with $i\Re(s)=\pm 1$.

\null Proof: \rm Suppose $\mu (\sigma\otimes\cdot )$ has a pole in
$\chi _{s\ti{\alpha }}$. Then, by results of Harish-Chandra
\cite{Si}, $\sigma $ is ramified, $P$ is self-conjugated,
$i_P^G\sigma $ is irreducible and $\mu(\sigma\otimes\chi _{\Im
(s)\ti{\alpha }})=0$. Write $\sigma _0=\sigma\otimes\chi _{\Im
(s)\ti{\alpha }}$. By \cite{Sh, 7.6}, there exists a unique
$i=1,2$ such that $0$ is a pole of $L(\cdot ,\sigma
_0,r_i)=L(\Im(s) +\cdot ,\sigma ,r_i)$. This proves the first
assertion and, by (1.3.3), $0$ is then also a pole of $L(\cdot
,\ti{\sigma}_0,r_i)$.

As $\chi _{\Re (s)\ti{\alpha }}$ is a pole of $\mu (\sigma
_0\otimes\cdot )$, it follows from the expression (1.3.2) for the
$\mu $-function and (1.3.5), that $1-i\Re(s)$ or $1+i\Re(s)$ is a
pole of $L(\cdot ,\ti{\sigma }_0,r_i)$, i.e. one of them must be
$0$. This concludes the proof of the first implication.

Conversely, choose $i$ such that $i\Re(s)=\pm 1$ and assume that
$\Im (s)$ is a pole of $L(\cdot ,\sigma ,r_i)$ (so that, in
particular, $i$ is an integer $\geq 1$). As $L(\cdot ,\sigma
,r_i)$ is regular for $i\geq 3$, we must have $i\in\{1,2\}$. Let
$\sigma _0=\sigma\otimes\chi _{\Im(s)\ti{\alpha }}$. Then $0$ is a
pole of $L(\cdot ,\sigma _0,r_i)$ by (1.3.4). By \cite{Sh, 7.4}
this can only happen if $P$ is self-conjugated. As the poles of
$L(\cdot ,\sigma _0,r_i)$ have real part $0$, it follows from
(1.3.2), (1.3.3) and (1.3.5) that $\mu (\sigma _0\otimes\cdot )$
has a pole in $\chi _{\Re(s)\ti{\alpha }}$.\hfill{\fin 2}

\null{\bf 2.2 Lemma:} \it Let $P=MU$ be a maximal standard
$F$-parabolic subgroup of $G$ and let $\psi _{\sigma
}:W_F\rightarrow\ ^LM$ be an elliptic admissible homomorphism.

Then, for any complex number $s$, the following two properties are
equivalent:

(i) $\Im (s)$ is a pole of $L(\cdot ,r_i\circ\psi _{\sigma })$ for
some positive integer $i$ verifying $i\Re(s)=\pm 1$;

(ii) $\ti{\alpha }(q)^{\Re(s)}$ is $q$-distinguished in the
connected centralizer of the image of the map $W_F\rightarrow\
^LM$, $\gamma\mapsto\ti{\alpha }(q)^{v_F(\gamma )\Im (s)}\psi
_{\sigma }(\gamma )$, and this connected centralizer is not a
torus.

\null Proof: \rm Replacing $\sigma $ by
$\sigma\otimes\chi_{\Im(s)\ti{\alpha }}$ and $\psi _{\sigma }$ by
$\gamma\mapsto\ti{\alpha }(q)^{v_F(\gamma )\Im (s)}\psi _{\sigma
}(\gamma )$, we can suppose by (1.3.4) that $\Im (s)=0$.

Denote by $\widehat{M}^{\sigma }$ the centralizer of $\psi
_{\sigma }(W_F)$ in $^LG$ and by $(\widehat{M}^{\sigma })^{\circ
}$ its connected component.

Suppose $\ti{\alpha }(q)^s$ $q$-distinguished in
$\widehat{M}^{\sigma }$ and that the connected component
$(\widehat{M}^{\sigma })^{\circ }$ is not a torus. So there is a
nilpotent element $N$ in the Lie algebra of the connected
component of $\widehat{M}^{\sigma }$, such that $(\Ad (\ti{\alpha
}(q)^s))(N)=qN$. Then $N\in V_{\pm i}$ for some integer $i$,
$1\leq i\leq m$, and it follows that $is=\pm 1$. Consider
$r_i\circ\psi _{\sigma }$. As $N\in V_{\pm i}^{W_F}$, we have
$V_{\pm i}^{I_F}\ne 0$. So the $L$-function $L(\cdot ,r_i\circ\psi
_{\sigma })$ is non trivial. As the Frobenius acts trivial on $N$,
it has a pole in $0$ by the above discussion of the Artin
$L$-function.

Conversely, choose $i$ such that $is=\pm 1$ and assume that $0$ is
a pole of $L(\cdot ,r_i\circ\psi _{\sigma})$ (so that in
particular $i$ is a positive integer). Replacing $s$ by $\vert
s\vert$, we can assume $is=1$. We will first prove that
$(\widehat{M}^{\sigma })^{\circ }$ is not a torus. As $L(\cdot
,r_i\circ\psi _{\sigma })$ is non trivial, there exists $N\in
V_i\subseteq Lie(^LU)$, which is invariant under the action of
$I_F$ by $r_i\circ\psi _{\sigma }$. As $L(\cdot ,r_i\circ\psi
_{\sigma })$ has a pole in $0$, by {\bf 1.4} we can choose $N$
such that the action of the Frobenius on $N$ by $r_i\circ\psi
_{\sigma }$ is trivial, i.e. $N$ is invariant by $W_F$. But then
$\exp (N)$ lies in the centralizer of $(r_i\circ\psi _{\sigma
})(W_F)$. So this centralizer contains a unipotent element. But,
the connected component of a reductive group which contains a
unipotent element cannot be a torus. So $(\widehat{M}^{\sigma
})^{\circ }$ is not a torus.

Remark that $(\Ad (\ti{\alpha }(q)^s))(N)=qN$. As
$\rk_{ss}(\widehat {M}^{\sigma })^{\circ }=1$, because
$T_{\widehat {M}}$, the maximal torus in the center of
$\widehat{M}$, is by \cite{H2, 4.2} a maximal torus of $(\widehat
{M}^{\sigma })^{\circ }$ and because $P$ (and consequently
$\widehat{M}$) is maximal, it follows that $\ti{\alpha }(q)^s$ is
$q$-distinguished in $(\widehat{M}^{\sigma })^{\circ
}$.\hfill{\fin 2}

\null {\bf 2.3 Theorem:} \it Let $G$ be the set of $F$-points of a
reductive connected quasi-split group defined over $F$, $P=MU$ a
maximal standard $F$-parabolic subgroup of $G$ and $\sigma $ a
unitary irreducible supercuspidal generic representation of $M$.

Let $\psi _{\sigma }:W_F\rightarrow\ ^LM$ be an admissible
elliptic homomorphism. Suppose that one has the equalities of
$L$-functions
$$L(\cdot,\sigma, r_i)=L(\cdot, r_i\circ\psi
_{\sigma }),\qquad 1\leq i\leq m.$$ Then the following property is
true:

For any complex number $s$, $\mu (\sigma\otimes\cdot )$ has a pole
in $\chi _{s\ti {\alpha }}$, if and only if $\ti {\alpha
}(q)^{\Re(s)}$ is $q$-distinguished in the connected centralizer
of the image of the map $W_F\rightarrow\ ^LM$,
$\gamma\mapsto\ti{\alpha }(q)^{v_F(\gamma )\Im (s)}\psi _{\sigma
}(\gamma )$, and this connected centralizer is not a torus.

Conversely, if this property is fulfilled and if all the poles of
the $L$-functions $L(\cdot ,r_i\circ\psi _{\sigma })$, $1\leq
i\leq m$, are simple, then the above equalities of $L$-functions
are true.

\null\it Proof: \rm If one has the equalities of $L$-functions,
the property in the theorem is a direct consequence of the lemmas
{\bf 2.1} and {\bf 2.2}. Conversely, if the property in the
theorem is true, the $L$-functions $L(\cdot ,\sigma ,r_i)$ and
$L(\cdot ,r_i\circ\psi _{\sigma })$ have the same poles on the
imaginary axes. So by (1.3.5) and {\bf 1.4} they have same poles
in $\Bbb C$. As $L(\cdot ,\sigma ,r_i)^{-1}$ and $L(\cdot
,r_i\circ\psi _{\sigma })^{-1}$ are both polynomials in $q^{-s}$
which take value $1$ in $0$, we conclude from the simplicity of
their zeroes (by (1.3.6) and by assumption) that they must be
equal.\hfill{\fin 2}

\null{\bf 2.4 Corollary:} \it In the notations and under the
assumptions of the preceding theorem assume that one has the
equality of $\gamma $-functions
$$\gamma (\cdot,\sigma, r_i,\psi _F)=\gamma (\cdot, r_i\circ\psi
_{\sigma },\psi _F),\qquad 1\leq i\leq m.$$

Then $\sigma $ verifies the assumption (LM) in \cite{H2, 4.3}
relative to $G$.\rm

\null\null{\bf 3.} In this section we will show that the
correspondence derived in \cite{H2} from the (conjectural) local
Langlands correspondence for supercuspidal representations
preserves $L$- and $\epsilon $-functions for generic
representations of quasi-split groups.

The following lemma is contained, but not explicitly stated in
\cite{Sh}.

\null{\bf 3.1 Lemma:} \it Let $G$ be the set of $F$-points of a
reductive connected quasi-split group, $P=MU$ and $P_1=M_1U_1$
standard $F$-parabolic subgroups of $G$, $M\supseteq M_1$, $\tau $
an irreducible smooth generic representation of $M_1$ and $\pi $
an irreducible smooth generic representation of $M_1$ which is a
sub-representation of $i_{P_1\cap M}^M\tau $.

Then, for any component $r_i$ of the adjoint representation $r:\
^LM\rightarrow Lie(^LU)$, we have
$$\gamma (\pi ,r_i,\psi _F)=\prod _{\alpha\in\sum _{red}(P_1),U_{1,\alpha
}\subseteq U}\gamma (\sigma ,r_{1,i,\alpha },\psi
_F),\eqno{(3.1.1)}$$ where $r_{1,i,\alpha }:\ ^LM_1\rightarrow
Lie(^LU_{1,\alpha })$ denotes the restriction of $r_i$.

\null\it Proof: \rm Denote by $\kappa _1$ the inclusion
$^LM_1\rightarrow\ ^LM$ and define $r_{1,i}=r_i\circ\kappa_1$. For
any root $\alpha\in\Sigma_{red}(P_1)$ verifying $U_{1 ,\alpha
}\subseteq U$, the space $Lie(^LU_{1 ,\alpha })$ is invariant by
$r_{1,i}$. Denote this representation by $r_{1,i,\alpha }$. Then
we have
$$r_{1,i}=\bigoplus _{\alpha\in\Sigma _{red}(P_1), U_{1
,\alpha }\subseteq U} r_{1,i,\alpha }.$$ Assume first $\tau $
supercuspidal. The product formula for the $\gamma $-function (cf.
[Sh, (3.13)]) gives an expression for $\gamma (\pi ,r_i,\psi _F)$
as a product of $\gamma $-functions related to $\tau $, which, by
the remarks in \cite{Sh, p. 306} (after the identity (6.2)) is in
fact a $\gamma $-factor attached to $\tau $ and $r_{1,i}$. The
unicity of that $\gamma $-factor and the identity \cite{Sh, p.
305} tell us that this $\gamma $-factor must be equal to $\prod
_{\alpha }\gamma (\tau ,r_{1,i,\alpha },\psi _F)$ with $\alpha\in
\Sigma _{red}(P_1)$, $U_{1,\alpha }\subseteq U$. The equality
(3.1.1) stated in the lemma follows.

If $\tau $ is no more supercuspidal, then there exist a standard
$F$-parabolic subgroup $P_2=M_2U_2$ of $G$, $M_2\subseteq M_1$,
and an irreducible supercuspidal representation $\sigma $ of
$M_1$, such that $\tau $ is a sub-representation of $i_{P_2\cap
M_1}^{M_1}\sigma $. By Theorem 2 of \cite{R}, $\sigma $ is
generic. Denote by $\kappa _2$ the inclusion $^LM_2\rightarrow\
^LM$, by $\kappa _{21}$ the inclusion $^LM_2\rightarrow\ ^LM_1$
and write $r_{2,i}=r_i\circ\kappa _2$ and
$r_{21,i}=r_{1,i}\circ\kappa _{21}$. Of course,
$r_{2,i}=r_{21,i}$. By, what we have just proved, we get that
$$\gamma (\pi ,r_i,\psi _F)=\prod _{\beta\in\Sigma
_{red}(P_2),U_{2,\beta }\subseteq U}\gamma (\sigma ,r_{2,i,\beta
},\psi _F)\eqno(*)$$ and, for $\alpha\in\Sigma_{red}(P_1)$, that
$$\gamma
(\tau ,r_{1,i,\alpha} ,\psi _F)=\prod _{\beta\in\Sigma
_{red}(P_2),U_{2,\beta }\subseteq U_{1,\alpha }}\gamma (\sigma
,r_{21,i,\beta },\psi _F).\eqno(**)$$ Substituting (**) in (*)
proves the identity (3.1.1) in the general case.\hfill{\fin 2}

\null {\bf 3.2 Theorem:} \it Let $G$ be the set of $F$-points of a
reductive connected quasi-split group defined over $F$. Let $\pi $
be a generic discrete series representation of a standard Levi
subgroup $M$ of $G$. Fix a standard parabolic subgroup
$P_1=M_1U_1$ of $G$ with $M_1\subseteq M$ and a generic
supercuspidal representation $\sigma $ of $M_1$ such that $\pi $
is a sub-representation of $i_{P_1\cap M}^M\sigma $.

Suppose that there is an admissible homomorphism $\psi _{\sigma
}:W_F\rightarrow\ ^LM_1$ such that for any $\alpha\in\Sigma
_{red}(P_1)$ and any irreducible component $r_{\alpha ,i}$ of
$r_{\alpha }:\ ^LM_{1,\alpha }\rightarrow\hbox{\rm Lie}(^LU_{1
,\alpha })$, we have
$$\gamma (r_{\alpha ,i}\circ\psi _{\sigma },\psi _F)=\gamma
(\sigma ,r_{\alpha ,i},\psi _F).$$ Then $\sigma $ verifies the
assumption $(LM)$ in \cite{H2}. Let $\psi _{\pi }$ be the
admissible homomorphisme $W_F\times\SL _2(\Bbb C)\rightarrow\ ^LM$
attached to $\pi $ in \cite{H2,5.3}. Then, for any component $r_i$
of the adjoint representation $r:\ ^LM\rightarrow Lie(^LU)$, we
have
$$\gamma (r_i\circ\psi _{\pi },\psi _F)=\gamma (\pi ,r_i, \psi
_F).$$

\null Proof: \rm It is a direct consequence of the corollary {\bf
2.4} that $\sigma $ verifies the assumption $(LM)$ in \cite{H2}
under our hypothesis. Denote by $P=MU$ the standard parabolic of
$G$ with Levi factor $M$ and by $r_{1,i}:\ ^LM_1\rightarrow
Lie(^LU)$ the restriction of $r_i$ to $^LM_1$. We have a
decomposition $r_{1,i}=\bigoplus _{\alpha\in\Sigma
(P_1),U_{1,\alpha}\subseteq U}r_{1,i,\alpha }$ with $r_{1,i,\alpha
}:\ ^LM_1\rightarrow Lie(^LU_{1,\alpha })$. Inserting our
assumptions $\gamma (\sigma ,r_{\alpha ,i},\psi _F)=\gamma
(r_{\alpha ,i}\circ\psi _{\sigma },\psi _F)$ in the identity
(3.1.1) and using the multiplicity of Artin $L$- and $\epsilon
$-functions, we get
$$\gamma (\pi ,r_i,\psi _F)=\gamma (r_{1,i}\circ\psi _{\sigma
},\psi _F).\eqno{(\#)}$$ Define $\psi _{\pi
}^{gal}:W_F\rightarrow\
^LM_1$ by $w\mapsto\psi _{\pi}(w,\pmatrix \vert w\vert^{1/2} & 0 \\
0 & \vert w\vert^{-1/2}\endpmatrix )$. (If one considers $\psi
_{\pi }$ as defined on the Weyl-Deligne group, then $\psi _{\pi
}^{gal}$ is the restriction of $\psi _{\pi }$ to $W_F$.) It is
proved in [Sh, 3.4] that $\gamma (r_i\circ\psi _{\pi },\psi
_F)=\gamma (r_i\circ\psi _{\pi }^{gal},\psi _F)$. As by
construction $\psi _{\pi }^{gal}=\psi _{\sigma }$, it follows that
$\gamma (r_{1 ,i}\circ\psi _{\sigma },\psi _F)=\gamma
(r_i\circ\psi _{\pi },\psi _F)$, which implies with the equality
(\#) the theorem. \hfill{\fin 2}

\null\null {\bf 3.3 Corollary:} \it Let $\pi $ be an irreducible
smooth generic representation of a standard Levi subgroup $M$ of
$G$. Fix a standard parabolic subgroup $P_1=M_1U_1$ of $G$ with
$M\supseteq M_1$ and an irreducible generic supercuspidal
representation $\sigma $ of $M_1$ such that $\pi $ is a
sub-representation of $i_{P_1\cap M}^{M}\sigma $.

Suppose that there is an admissible homomorphism $\psi _{\sigma
}:W_F\rightarrow\ ^LM_1$ such that for any $\alpha\in\Sigma
_{red}(P_1)$ and any irreducible component $r_{\alpha ,i}$ of
$r_{\alpha }:\ ^LM_{1,\alpha }\rightarrow\hbox{\rm Lie}(^LU_{1
,\alpha })$, we have
$$\gamma (r_{\alpha ,i}\circ\psi _{\sigma },\psi _F)=\gamma
(\sigma ,r_{\alpha ,i},\psi _F).$$ Then $\sigma $ verifies the
assumption $(LM)$ in \cite{H2}. Let $\psi _{\pi }$ be the
admissible homomorphism $W_F\times\SL _2(\Bbb C)\rightarrow\ ^LM$
attached to $\pi $ in \cite{H2, 5.5}. Then, for any component
$r_i$ of the adjoint representation $r:\ ^LM\rightarrow Lie(^LU)$,
we have
$$L (\cdot ,r_i\circ\psi _{\pi })=L (\cdot ,\pi ,r_i)\qquad
\hbox{\it and}\qquad \epsilon(r_i\circ\psi _{\pi },\psi
_F)=\epsilon(\pi ,r_i, \psi _F).$$

\null\it Proof: \rm We will first consider the case, when $\pi $
is tempered. Then it is by {\bf 1.5} enough to show that
$$\gamma(r_i\circ\psi _{\pi },\psi_F)=\gamma(\pi ,r_i, \psi
_F)$$ for any $i$. After possibly changing $\sigma $ (and
consequently $\psi _{\sigma }$) by an unramified character twist
(which conserves by (1.3.4) and {\bf 1.4} the equalities of
$\gamma $-functions), we can find a standard parabolic subgroup
$P_2=M_2U_2$, $M\supseteq M_2\supseteq M_1$, and a generic
irreducible discrete series representation $\tau $ of $M_2$ which
is a sub-representation of $i_{M_2\cap P_1}^{M_2}\sigma $, such
that $\pi $ is a sub-representation of $i_{M\cap P_2}^M\tau $.
Denote by $\kappa _2$ the inclusion $^LM_2\rightarrow\ ^LM$ and
put $r_{2,i}=r_i\circ \kappa _2$. By the identity (3.1.1) and
Theorem {\bf 3.2} we have
$$\eqalign{\gamma (\pi ,r_i,\psi _F)=&\prod _{\alpha\in\Sigma_{red}(P_2),
U_{2,\alpha }\subseteq U}\gamma (\tau ,r_{2,i,\alpha },\psi _F)\cr
=&\prod _{\alpha\in\Sigma _{red}(P_2),U_{2,\alpha }\subseteq
U}\gamma (r_{2,i,\alpha }\circ\psi _{\tau },\psi _F)\cr =&\gamma
(r_{2,i}\circ\psi _{\tau },\psi _F).\cr}$$ As by construction
$\psi _{\tau }$ and $\psi _{\pi }$ take the same values, it
follows that $\gamma(r_i\circ \psi _{\pi },\psi_F)=\gamma(\pi
,r_i, \psi _F)$.

Let now $\pi $ be an arbitrary generic smooth representation of
$M$. Then, after possibly changing $\sigma $ and $\psi _{\sigma }$
by an unramified character twist, using Langlands' classification,
there is a semi-standard parabolic subgroup $P_2=M_2U_2$ of $G$
with $M\supseteq M_2\supseteq M_1$ and a generic quasi-tempered
representation $\tau $ of $M_2$ such that $\tau $ is a
sub-representation of $i_{P_1\cap M_2}^{M_2}\sigma $ and $\pi $ is
a sub-representation of $i_{P_2\cap M}^M\tau $. By (1.3.7), $L
(\cdot ,\pi ,r_i)$ is a product of $L$-functions attached to $\tau
$ with respect to simple reflections of $P_2$. As $L(\cdot
,r_i\circ \psi _{\pi })$ is obtained in the same way from the
$L$-functions of $\psi _{\tau }$, the equality of the
$L$-functions of $\pi $ and $\psi _{\pi }$ follows from the
tempered case proved just before. The proof of the equality of
$\gamma $-functions is literally the same as for $\pi $ tempered.
The identity for $\epsilon $-factors follows from this (cf. {\bf
1.3} and {\bf 1.5}).\hfill{\fin 2}

\null\null{\bf 4.} We will now finish with remarks on the general
case, i.e. we will consider representations which are not generic
and later also groups which are not quasi-split.

\null{\bf 4.1} So suppose first that $G$ is still the set of
$F$-points of a quasi-split connected reductive group. In order to
define $L$-functions and $\epsilon $-factors for non generic
representations, two assumptions are made in \cite{Sh} (and
justified by other more basic assumptions).

\null (4.1.1) Each tempered $L$-packet of a standard Levi subgroup
contains a generic representation.

\null (4.1.2) Harish-Chandra's $\mu $-function defined on discrete
series depends only on $L$-packets.

\null Let $P=MU$ be a standard $F$-parabolic subgroup of $G$ and
$\pi $ a non generic irreducible tempered representation of $M$.
Let $r_i$ be a component of the adjoint representation $r:\
^LM\rightarrow Lie(^LU)$. By assumption (4.1.1) there exists a
generic irreducible representation $\pi '$ in the $L$-packet of
$\pi $. One defines $L(\cdot ,\pi ,r_i)$, $\epsilon (\cdot ,\pi
,r_i,\psi _F)$ and $\gamma (\cdot ,\pi ,r_i,\psi _F)$ to be
$L(\cdot ,\pi ',r_i)$, $\epsilon (\cdot ,\pi ',r_i,\psi _F)$ and
$\gamma (\cdot ,\pi ',r_i,\psi _F)$ respectively.

Let now $\pi $ be an arbitrary irreducible smooth non generic
representation of $M$. By Langlands' classification, there is a
standard $F$-parabolic subgroup $P_1=M_1U_1$ of $G$, $M\supseteq
M_1$, and an irreducible quasi-tempered representation $\tau $ of
$M_1$, such that $\pi $ is the unique sub-representation of
$i_{P_1\cap M}^M\tau $. By assumption (4.1.1) there is a generic
quasi-tempered representation $\tau '$ in the $L$-packet of $\tau
$ such that $i_{P_1\cap M}^M\tau '$ has a unique
sub-representation $\pi '$. This representation $\pi '$ may not be
generic, but we define $L(\cdot ,\pi ',r_i)$, $\epsilon (\cdot
,\pi ',r_i,\psi _F)$ and $\gamma (\cdot ,\pi ',r_i,\psi _F)$ by
the same formulas (see (1.3.7) and following) as in the generic
case. The local factors for $\pi $ are by definition those for
$\pi '$.

\null {\bf 4.2} To extend our discussion of the results in
\cite{H2} to non generic representations of $G$, we have in order
to use the results in section {\bf 2} and {\bf 3} to make the
following assumption.

\null (4.2.1) If $\sigma $ is an irreducible generic supercuspidal
representation of an $F$-Levi subgroup $M$ of $G$, then there is
an admissible homomorphism $\psi _{\sigma }:W_F\rightarrow\ ^LM$,
which has the same local factors with respect to the adjoint
action of $^LM$ then $\sigma $.

\null {\bf 4.3} We are now able to deduce from this the general
version of the assumption (PM) in \cite{H2}:

\null {\bf Lemma:} \it Let $G$ be the set of $F$-points of a
quasi-split group, $P=MU$ a standard parabolic subgroup of $G$ and
$\sigma $ a unitary irreducible supercuspidal representation of
$M$. Suppose that (4.1.1) and (4.2.1) hold. Then there is a
discrete series representation $\tau $ in the $L$-packet of
$\sigma $, a standard parabolic subgroup $P_1=M_1U_1$ of $G$,
$M\supseteq M_1$, and an irreducible unitary supercuspidal
representation $\sigma _1$ of $M_1$, $\tau \subseteq i_{P_1\cap
M}^M(\sigma _1\otimes\chi _{\lambda })$ for some $\lambda \in
a_{M_1}^*$, with the following property with respect to any root
$\alpha\in\Sigma(P_1)$:

Let $s$ be a real number $s>0$. Then $\mu (\sigma_1\otimes\cdot )$
has a pole in $\chi _{s\ti {\alpha }}$, if and only if $\ti
{\alpha }(q)^s$ is $q$-distinguished in the connected centralizer
of $\psi _{\sigma }(W_F)$ and this connected centralizer is not a
torus.

In addition, one can choose for $\sigma _1$ a generic
representation.

\null\it Proof: \rm By (4.1.1), there is a generic representation
$\tau $ in the $L$-packet of $\sigma $, which must be a discrete
series. One can choose $P_1=M_1U_1$ as in the statement and an
irreducible supercuspidal representation $\sigma _1$ of $M_1$ such
that $\tau $ is a sub-representation of $i_{P_1\cap M}^M\sigma
_1$. The representation $\sigma _1$ must be generic by \cite{R,
Theorem 2}. So, using the assumption (4.2.1), the corollary {\bf
2.4} applies and proves the theorem. \hfill{\fin 2}

\null {\bf 4.4 Theorem:} \it Let $G$ be the set of $F$-points of a
quasi-split group, $P=MU$ a standard parabolic subgroup of $G$ and
$\pi $ an irreducible smooth representation of $M$. If the
assumptions (4.1.1), (4.1.2) and (4.2.1) are verified, then the
construction in \cite {H, 5.5}, that associates to $\pi $ an
admissible homomorphism $\psi _{\pi }:W_F\times\SL_2(\Bbb
C)\rightarrow\ ^LM$, applies and one has
$$L (\psi _{\pi }\circ r_i)=L (\pi ,r_i)\qquad\hbox{\rm and}\qquad
\epsilon(\psi _{\pi }\circ r_i,\psi _F)=\epsilon(\pi ,r_i, \psi
_F).$$

\null\it Proof: \rm By Langlands' classification, there is a
standard $F$-parabolic subgroup $P_1=M_1U_1$ of $G$, $M\supseteq
M_1$, and an irreducible quasi-tempered representation $\tau $ of
$M_1$, such that $\pi $ is the unique sub-representation of
$i_{P_1\cap M}^M\tau $. By assumption (4.1.1) there is a generic
representation $\tau '$ in the $L$-packet of $\tau $. The unique
sub-representation $\pi '$ of $i_{P_1\cap M}^M\tau '$ is in the
same $L$-packet than $\pi $. By {\bf 4.1} it has the same local
factors than $\pi $. The proof of corollary {\bf 3.3} generalizes
to $\pi '$, showing that the local factors for $\pi '$ and $\psi
_{\pi '}$ are the same. As one can choose $\psi _{\pi }=\psi _{\pi
'}$, this proves the theorem.\hfill{\fin 2}

\null {\bf 4.5} Consider now that $G$ is the set of $F$-points of
an arbitrary connected reductive group defined over $F$ which may
not be quasi-split. It is believed that Harish-Chandra's $\mu
$-function is invariant for inner forms (cf. \cite{Sh, 9}). The
constructions in \cite{H2} are also invariant for inner forms.
Local factors for representations of Levi subgroups of $G$ are
defined by the ones for the corresponding representations for the
quasi-split inner form of $G$. So it is clear that the
correspondence must conserve the local factors.

\null \sit Acknowledgement: \srm The author started to work on the
subject of the article after a one year-stay at Purdue University.
He took some profit from discussions with F. Shahidi during that
time. The author is in particular indebted to F. Shahidi for
having had a look on a final version of this paper. Thanks also to
A.-M. Aubert for pointing out some misprints.\rm

\bye